\documentclass{svjour3_VERWURSCHTET}
\usepackage{graphicx,xcolor}
\newcommand{\ii}{{\mathrm i}}
\def\N     {\ensuremath{\mathbf N}}
\def\R     {\ensuremath{\mathbf R}}

\newcommand{\ee}{{\mathrm e}}

\newcommand{\nE}{\mathcal{E}}
\newcommand\EA{\nE_A}
\newcommand\EB{\nE_B}

\newcommand\rev[1]{{{#1}}}

\begin{document}
\title{Adaptive Exponential Integrators for MCTDHF\thanks{Supported by the Vienna Science and Technology Fund (WWTF) grant MA14-002.
The computations have been conducted on the Vienna Scientific Cluster (VSC).}}
\author{Winfried Auzinger, Alexander Grosz,
Harald Hofst{\"a}tter, and Othmar Koch
}
\authorrunning{W. Auzinger et al.}
\institute{Vienna University of Technology, Wiedner Hauptstra{\ss}e 8--10, A-1040 Wien, Austria
%\email{w.auzinger@tuwien.ac.at, e1525490@student.tuwien.ac.at} \and
%University of Vienna, Institute of Mathematics, Oskar-Morgenstern-Platz 1, A-1090 Wien, Austria
%\email{hofi@harald-hofstaetter.at, othmar@othmar-koch.org}
}
\maketitle              % typeset the header of the contribution
\begin{abstract}
We compare exponential-type integrators for the numerical time-propagation of the
equations of motion arising in the multi-configuration time-dependent Hartree-Fock
method for the approximation of the high-dimensional multi-particle Schr{\"o}dinger equation.
We find that among the most widely used integrators
like Runge-Kutta,
exponential splitting, exponential Runge-Kutta,
exponential multistep and Lawson methods, exponential
Lawson multistep methods with one predictor/corrector step provide optimal
stability and accuracy at the least computational cost, taking into account that
the evaluation of the nonlocal potential terms is by far the computationally most
expensive part of such a calculation. Moreover, the predictor step provides an
estimator for the time-stepping error at no additional cost, which enables
adaptive time-stepping to reliably control the accuracy of a computation.

\end{abstract}

\section{Introduction} \label{sec:intro}

We compare  time integration methods for nonlinear Schr{\"o}dinger-type
equations %of the type
\begin{equation}\label{de1}
\ii\,\partial_t u(t) =  A\,u(t) + B(u(t)) = H(u(t)),~\,t>t_0,
\quad u(t_0)=u_0,
\end{equation}
on \rev{the Hilbert} space $\mathcal{B}\rev{=L^2}$.
Here, $A\colon \mathcal{D}\subseteq \mathcal{B} \to \mathcal{B}$ is a self-adjoint
differential operator and $B$ a generally unbounded nonlinear operator.
Our focus is on the equations of motion associated with the multi-configuration
time-dependent Hartree-Fock
(MCTDHF) approximation to the multi-particle electronic Schr{\"o}dinger equation,
where the key issue is the high computational effort for the evaluation of the
nonlocal (integral) operator $B$. Thus, in the choice of the most appropriate integrator,
we emphasize a minimal number of evaluations of $B$ for a given order and
disregard the effort for the propagation of $A$, which can commonly be realized
at essentially the cost of two (cheap) transforms between real and frequency space
via fast transforms like [I]FFT.
%It turns out that splitting methods require a prohibitive number of evaluations
%of $B$ to obtain a high-order approximation, whence in spite of their favorable properties, splitting methods are computationally too expensive for our purpose.
The approaches that we pursue and advocate in this paper are thus based on splitting of the vector fields in~(\ref{de1}).
\rev{It turns out that exponential integrators \cite{hocost10} based on the variation of constants serve our purpose best, as they
provide a desirable balance between computational effort and stability.}

\section{The MCTDHF method}
\label{subsec:mctdhf}

We focus on the comparison of numerical methods for the equations of
motion associated with MCTDHF for the approximate solution of the time-dependent
multi-particle Schr{\"o}dinger equation
$$ %\begin{equation}\label{eq:schroe1}
\ii\,\frac{\partial \psi}{\partial t}=H\psi,
$$ %\end{equation}
where the complex-valued wave function
$\psi=\psi(x_1,\dots,x_f,t)$ depends on time~$t$ and,
in the case considered here, the positions $x_1,\dots,x_f\in\R^3$ of
electrons in an atom or molecule. The time-dependent Hamiltonian reads
\begin{eqnarray*}%\label{eq:hamilton1}
H&=&H(t):=\sum_{k=1}^{f}\Big(\,\frac12\,\Delta^{(k)} +
U(x_k)+\sum_{\ell<k} V(x_k-x_\ell)\Big)+ V_{\mathrm{ext}}(x_1,\dots,x_f,t)\nonumber\\
&=:& T + W(t,x_1,\dots,x_f),\\
T&=&\sum_{k=1}^{f}\frac12\,\Delta^{(k)},
\quad U(x)=-\frac{Z}{|x|}, %=-\frac{Z}{\sqrt{x_1^2+x_2^2+x_3^2}},
%\label{eq:potential1}
~Z\in \N,
\quad V(x-y)=\frac{1}{|x-y|}\,.
%=\frac{1}{\sqrt{(x_1-y_1)^2+(x_2-y_2)^2+(x_3-y_3)^2}},
%\label{eq:potential2}
\end{eqnarray*}
Here $V_{\mathrm{ext}}(x_1,\dots,x_f,t)$ is a smooth time-dependent function, and $\Delta^{(k)}$ is the Laplace operator
with respect to $x_k$ only.

In MCTDHF as put forward in~\cite{zanghellinietal03c},
the multi-electron wave function $\psi$
%from~(\ref{eq:schroe1})
is approximated by a function $u$ living in a manifold $\mathcal{M}$
characterized by the ansatz
\begin{equation}\label{eq:ansatz1}
u = \sum_{(j_1,\dots,j_f)} a_{j_1,\dots,j_f}(t)\,\phi_{j_1}(x_1,t)
\cdots\phi_{j_f}(x_f,t)=:\sum_{J}a_J(t)\,\Phi_J(x,t).
\end{equation}
%
%Using~(\ref{eq:ansatz1})
For the electronic Schr\"odinger
equation, the Pauli principle implies that
only solutions $u$ are considered
which are antisymmetric under exchange of any pair of
arguments $x_j, x_k$,
%This assumption is particular to the MCTDHF approach,
%as compared to the multi-configuration time-dependent Hartree
%method (MCTDH) proposed in~\cite{becketal00,becmey97,meyeretal90,meyerworth03} for quantum molecular dynamics.
%reducing the number of equations for $a_J$ to ${N\choose f}$.

%In the case of the electronic Schr\"odinger equation which we are
%focussing on here, we further have to take into account electron
%spin. However, it was explained for example in~\cite{kochetal04a}
%that this does not change our considerations concerning the equations
%of motion associated with MCTDHF.
%Consequently, we ignore spin and concentrate on the representation
%(\ref{eq:ansatz1}).

Now, the Dirac-Frenkel variational principle
\cite{frenkel34} \rev{in conjunction with orthogonality conditions} is used to derive differential
equations for the coefficients $a_J$ and the so-called single-particle functions
$\phi_{j}$ in~(\ref{eq:ansatz1}), \rev{where we will henceforth tacitly identify $u$ with
the vector $(a,\phi)$ of coefficients and orbitals,}
%This means that for $u\in\mathcal{M}$ we require
%%
%\begin{equation}\label{eq:diracfrenkel1}
%\big\langle \delta u\big|\,\ii\textstyle{\frac{\partial}{\partial t}}
%-H\big|u\big\rangle=0,
%\end{equation}
%%
%where $\delta u$ varies in the tangent space
%${\mathcal T}_u\,{\mathcal M}$ of ${\mathcal M}$ at $u$.
%%
%%It was shown in~\cite{koclub05a} that the set $\mathcal{M}$
%%in conjunction with a full-rank-condition for the
%%\emph{density matrix} $\rho$ defined in~(\ref{eq:density}) below can
%%be endowed with the structure of a manifold, justifying the
%%application of the variational principle as explained above.
%%We do not give details here, but henceforth sloppily
%%refer to $\mathcal{M}$ as a manifold under the assumption that
%%$\rho$ is nonsingular.
%%
%In order to define a unique solution of~(\ref{eq:diracfrenkel1})
%additional constraints are imposed,
%%
%\begin{equation}
%\langle \phi_j\big|\phi_k\rangle =
%\delta_{j,k},~\, t\ge 0,\quad
%\big\langle \phi_j\big| \textstyle{\frac{\partial \phi_k}{\partial t}}\big\rangle
%= -\ii\,\langle \phi_j\,|\,T\,|\,\phi_k\big\rangle.\label{eq:ortho2}
%\end{equation}
%%
%The variational principle~(\ref{eq:diracfrenkel1}) and the additional
%restrictions~(\ref{eq:ortho2}) %finally
%yield equations of motion for the coefficients $ a_J $ and the
%single-particle functions $ \Phi_J $ in~(\ref{eq:ansatz1}),
%
\begin{eqnarray}
&& \ii\,\frac{da_J}{dt}
=\sum_K\big\langle \Phi_J\left| W \right|
\Phi_K\big\rangle\,a_K \quad \forall J,\label{eq:work1}\\
&& \ii\,\frac{\partial \phi_j}{\partial t}= T\,\phi_j+
(I-P)\sum_{k=1}^N\sum_{\ell=1}^N\rho^{-1}_{j,\ell}\,\overline{W}_{\ell,k}\,\phi_k,\quad
j=1,\dots,N,\label{eq:work2}
\end{eqnarray}
where
$$ %\label{eq:meanfield}
\overline{W}_{j,\ell} = \big\langle \psi_j|W|\,\psi_\ell\big\rangle,
\quad \mbox{with} \quad
\psi_j=\big\langle \phi_j\,|\,u\big\rangle,~~
\rho_{j,\ell} = \big\langle \psi_j\,|\,\psi_\ell\big\rangle,
$$
and where $P$ is the orthogonal projector onto the space spanned by the
functions~$\phi_{j}$.
%Clearly, the high-dimensional integrals appearing
%in~(\ref{eq:work1})--(\ref{eq:work2}) imply a huge computational effort
%for each evaluation, whence the number of evaluations of $B$ should be
%kept to a minimum as detailed above.
%
We will henceforth denote
\begin{equation}
A=\frac{\ii}{2}\Big(0,\Delta^{(1)},\dots,0,\Delta^{(f)}\Big)^T,
\quad B=B(a,\phi),
\end{equation}
where $B$ is the vector of the components associated with the potential %to correspond with (\ref{de1}),
%and we denote the coefficient tensor by $a=(a_J)$ and the vector of orbitals by $\phi=(\phi_j)$.
\rev{which constitute the computationally most expensive part.}

\subsection{Splitting methods} Popular integrators for quantum dynamics are exponential
time-splitting methods which are based on multiplicative combinations
of the partial flows $ \EA(t,u)\colon$ $ u \mapsto u(t)=\ee^{tA}\,u $ and
$ \EB(t,u)\colon u \mapsto u(t) $ with $u'(t)=B(u(t)),\;u(0)=u$.
For a single step $ (t_n,u_n) \mapsto (t_{n+1},u_{n+1}) $ with time-step $h $, this reads
$$%\label{splitting1}
u_{n+1} := \mathcal{S}(h,u_n) = \EB(b_s h,\cdot) \circ \EA(a_s h,\cdot) \circ \ldots \circ
                              \EB(b_1 h,\cdot) \circ \EA(a_1 h,u_n),
$$
where the coefficients $ a_j,b_j,\, j=1 \ldots s $
are determined according to the requirement
that a prescribed order of consistency is obtained~\cite{haireretal02}.
%
%A rigorous mathematical error analysis of splitting methods was
%first given in~\cite{jahlub00} for low order methods applied to linear Schr\"{o}dinger
%equations. An extension to higher-order schemes is provided by~\cite{th12}. The nonlinear
%setting was first analyzed for the second-order Strang splitting in~\cite{lubich07},
%see also~\cite{koclub08c}, and a convergence proof for high-order
%methods which also covers the cubic nonlinear Schr{\"o}dinger equation and the equations
%of motion associated with MCTDHF is given in~\cite{knth10a}.
%Theoretical error bounds for rotating BECs by a Fourier--Laguerre--Hermite splitting method
%are given in~\cite{hofkoctha12}, while fully implicit finite difference methods are analyzed in~\cite{baocai12}.
%An alternative theoretical framework for the analysis of splitting methods based on the defect of the numerical solution was
%recently developed in~\cite{auzingeretal13a} for linear problems and also extended to the nonlinear case
%in~\cite{auzingeretal13b}.
%The convergence was analyzed for evolution equations of Schr{\"o}dinger type.
%This also includes the construction and theoretical analysis of a~posteriori error estimators
%for the purpose of designing adaptive schemes.
%
For a convergence analysis of splitting methods in the context of MCTDHF, see for instance~\rev{\cite{knth10a}}.

\subsection{Exponential integrators}
An approach which also exploits the separated vector fields is given by
the class of exponential integrators, which %~\cite{hersch58,hochbrucketal98}.
are comprehensively discussed in~\cite{hocost10}.
Here the variation of constant formula is used to express the solution
of~(\ref{de1}) for a time-step $t_n \to t_{n+1}=t_n+h$ via the integral equation
\begin{equation}\label{voc}
u(t_n+h) =
\ee^{hA}\,u_n + \int_0^h\ee^{(h-\tau)A}B(u(t_n+\tau))\,\mathrm{d}\tau.
\end{equation}
Different numerical integrators are distinguished depending on how the
integral in~(\ref{voc}) is approximated.

\paragraph{Exponential Runge-Kutta methods}

When the integral in~(\ref{voc}) is approximated by a quadrature formula of Runge-Kutta type,
relying on evaluations of the nonlinear operator $B$ at interior points $t_n + h\tau_j,\
\tau_j\in [0,1],\ j=1,\dots,k$,
an \emph{exponential Runge-Kutta} method is obtained.
This corresponds to replacing $B(\cdot)$ in the integrand
by a polynomial interpolant at the points
$$
\big( t_n+h\tau_1,B(u(t_n+h\tau_1) \big),\dots,
\big( (t_n+h\tau_k,B(u(t_n+h\tau_k) \big).
$$
The method is realized by stepping from $t_n+h\tau_j \to t_n+h\tau_{j+1}$ in the same way as for
a Runge-Kutta method, with appropriate weights of the stages. For implicit methods, nonlinear
systems of equations have to be solved, which is generally considered as prohibitive.
Note that after interpolation, the resulting integral can be evaluated analytically by using the $\varphi$-functions or alternatively, by numerical quadrature~\cite{hocost10}.
%see for example~\cite{niewri12}.
Such a procedure has first been proposed in~\cite{friedli78},
for a stiff error analysis, see~\cite{hocost10} and references therein.
For our comparisons, we use the fourth order Krogstad method mentioned there.
%The denotation \emph{exponential integrator} first appeared in~\cite{hochbrucketal98},
%where explicit exponential Runge-Kutta methods are considered. The non-stiff order
%conditions are derived, and fourth order methods are constructed along with reduced methods
%which promise significant computational advantages when used in conjunction with
%Krylov methods for the matrix exponential. Also, stepsize control linked with the Krylov-substeps is proposed,
%and embedded Runge-Kutta methods serve as a basis for adjusting the time step. More recently,
%a systematic theory to derive the order conditions based on the variation of constant
%formula and trees similarly as in classical Runge-Kutta methods has been
%established in~\cite{luanost13}. The latest developments to our knowledge are studies
%of exponential Crank-Nicolson integrators in conjunction with Pad{\'e} approximation
%of the matrix exponential given in~\cite{liangetal14}, which also gives a linear stability
%analysis and uses extrapolation to improve the basic order two, a comparative study
%of exponential integrators~\cite{monboo17} which seems to favor explicit exponential Runge-Kutta
%methods, and a similar investigation for a cubic--quintic complex Ginzburg--Landau equation,
%where explicit Lawson Runge-Kutta methods show the best performance as compared
%to explicit Runge-Kutta or splitting methods. An adaptive approach based on embedded
%pairs of formulae is also studied there.

\paragraph{Exponential multistep methods}

The integral in~(\ref{voc}) can be approximated in terms of an interpolation polynomial at
previous approximations
\begin{equation}\label{twostar}
\big( -(k-1)\,h,B(u_{n-k+1}) \big),\dots,
\big((-h,B(u_{n-1})),(0,B(u_n) \big)\,.
\end{equation}
This yields an (explicit) exponential Adams-Bashforth multistep method \rev{first mentioned in \cite{certaine60},
and introduced more systematically  in \cite{noersett69}, see also for instance~\cite{calvopal11,hocost10}}.
If the interpolation also comprises the forward point $(h, B(u_{n+1}))$,
%
%$$(-(k-1)h,B(u_{n-k+1})),\dots, (-h, B(u_{n-1})),(0, B(u_n)), (h, B(u_{n+1})),$$
%
an (implicit) exponential Adams-Moulton method is obtained. These two approaches can be combined in
a predictor/corrector method in the same way as for linear multistep methods.
Exponential multistep methods have first been considered and analyzed in~\cite{calvopal11} under
the assumption of smooth $B$, and a starting strategy is also given there.
%A similar approach is adopted in~\cite{hocost11},
%with a different integration domain in the variation of constant formula, and
%a linearized version is proposed. The stability analysis proceeds as in~\cite{hocostsch09}.
%Related exponential integrators of a multistep flavor are also considered in
%\cite{ostermannetal06,auzlap12}.

\paragraph{Lawson methods}

In Lawson methods, equation~(\ref{de1}) is transformed prior
to the numerical integration by the substitution
$ u(t) \to \ee^{-tA}\,u(t)$.
To the resulting equation
\begin{equation}\label{lsmf1}
u'(t) = \ee^{-tA}B\,(\ee^{tA}\,u(t)) =: F(u(t)),
\end{equation}
any appropriate time-stepping scheme can be
applied.
%\footnote{%
%In our numerical experiments reported below,
%we have solved the transformed problem by the classical explicit
%fourth order Runge-Kutta method (RK4).}
%\rev{In a one-step version, an explicit Runge--Kutta method is employed to solve (\ref{lsmf1}),
%which is equivalent to interpolation at interior nodes of the whole integrand
%in (\ref{voc}) by a polynomial. Lawson multistep methods result from interpolation of the whole integrand
%in (\ref{voc}) at previous points in the same fashion as in (\ref{twostar}).}%
%}
The main advantage lies in
the fact that the dynamics associated with the non-smooth operator $A$ is separated
by the transformation which can be realized cheaply in frequency space, while the
problem subjected to the time-stepping scheme is smoother, thus allowing for
larger time-steps.
This transformation was first introduced in~\cite{lawson67}
for ordinary differential equations.

\rev{In a one-step version, an explicit Runge--Kutta method is employed to solve
(\ref{lsmf1}), which is equivalent to interpolation at interior nodes of the
whole integrand in (\ref{voc}) by a polynomial in the same fashion as in (\ref{twostar}).}
Reference~\cite{hocost17} gives a convergence proof of Lawson Runge-Kutta methods in the
stiff case, however under the assumption that the operator $B$ is smooth, which
is not the case in the MCTDHF equations we are considering. A convergence proof for Adams-Lawson multistep
methods for the MCTDHF equations under minimal regularity requirements is given in the forthcoming work~\cite{koch19}.
\rev{The proof addresses the transformed equation (\ref{lsmf1}) and combines stability
and consistency to conclude convergence. To this end, a boot-strapping
argument is employed, first showing convergence in the Sobolev space $H^1$.
Stability in $L^2$ only holds if the numerical solution is in $H^1$, which
follows from the first argument, whence convergence
in $L^2$ is inferred. Lipschitz conditions for the right-hand side entering the stability arguments
can be shown by appropriate Sobolev-type inequalities in both $H^1$ and $L^2$.
To prove consistency, the norms of derivatives of $F$ in (\ref{lsmf1}) are estimated,
which amounts to bounds on commutators of the operators $A$ and $B$. This implies
assumptions on the regularity of the exact solution $u$.}

%Alternatively, we can view the Lawson approach as an approximation to~(\ref{voc}),
%where the whole integrand is replaced by a polynomial interpolant. We are mainly
%interested in multistep methods, as evaluations of $B$ are prohibitively expensive,
%so exponential Lawson multistep are constructed by replacing the integrand by an interpolant at
%$$
%(-(k-1)h,\ee^{(k-1)hA}B(u_{n-k+1})),\dots, (-h, \ee^{hA}B(u_{n-1})),(0, B(u_n)),\dots, [(h,\ee^{-hA}B(u_{n+1}))].
%$$
%The last point in square brackets is used for implicit methods, while otherwise
%the methods are explicit.
We will demonstrate that the best approach for our
goal is to use exponential Lawson multistep methods in a
predictor/corrector implementation,
which is shown to increase the accuracy
%of the numerical method
and also provides a local error
estimator for adaptive time-stepping at no additional cost.
The efficiency of the time discretization can be improved if high-order time propagators are employed.
In the multistep approach, this does not imply additional computational cost
if no memory limitations have to be taken into account.

\paragraph{Comparisons}
To assess the performance of the exponential integration methods described above, we will
also show results for the classical explicit Runge-Kutta method of fourth order (RK4)
\rev{and the second-order Strang splitting}.

\section{Numerical results}\label{sec:num}

\begin{figure}[ht!]
\begin{center}
\includegraphics[width=0.482\textwidth]{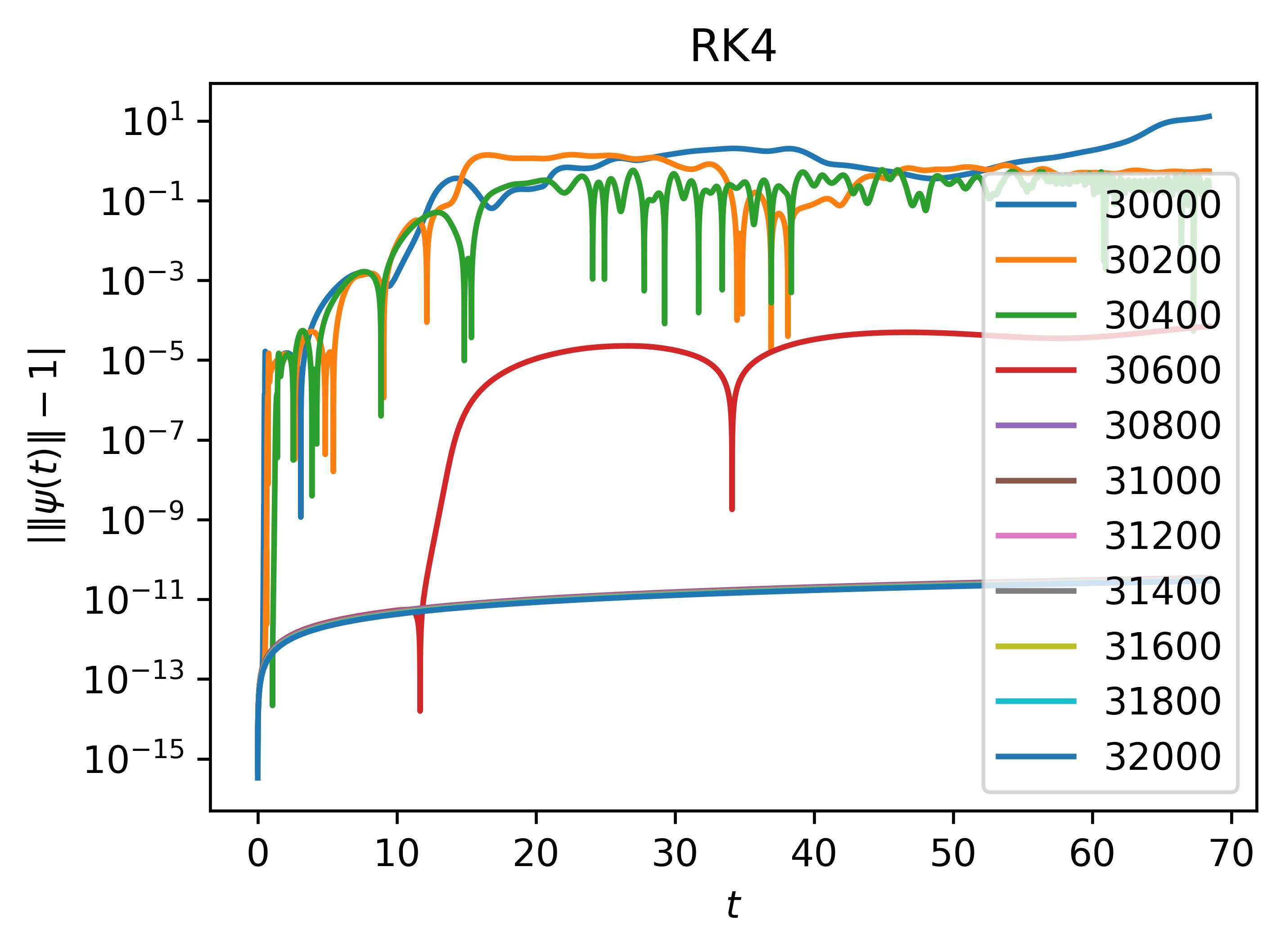}
\hspace{0.02\textwidth}
\includegraphics[width=0.482\textwidth]{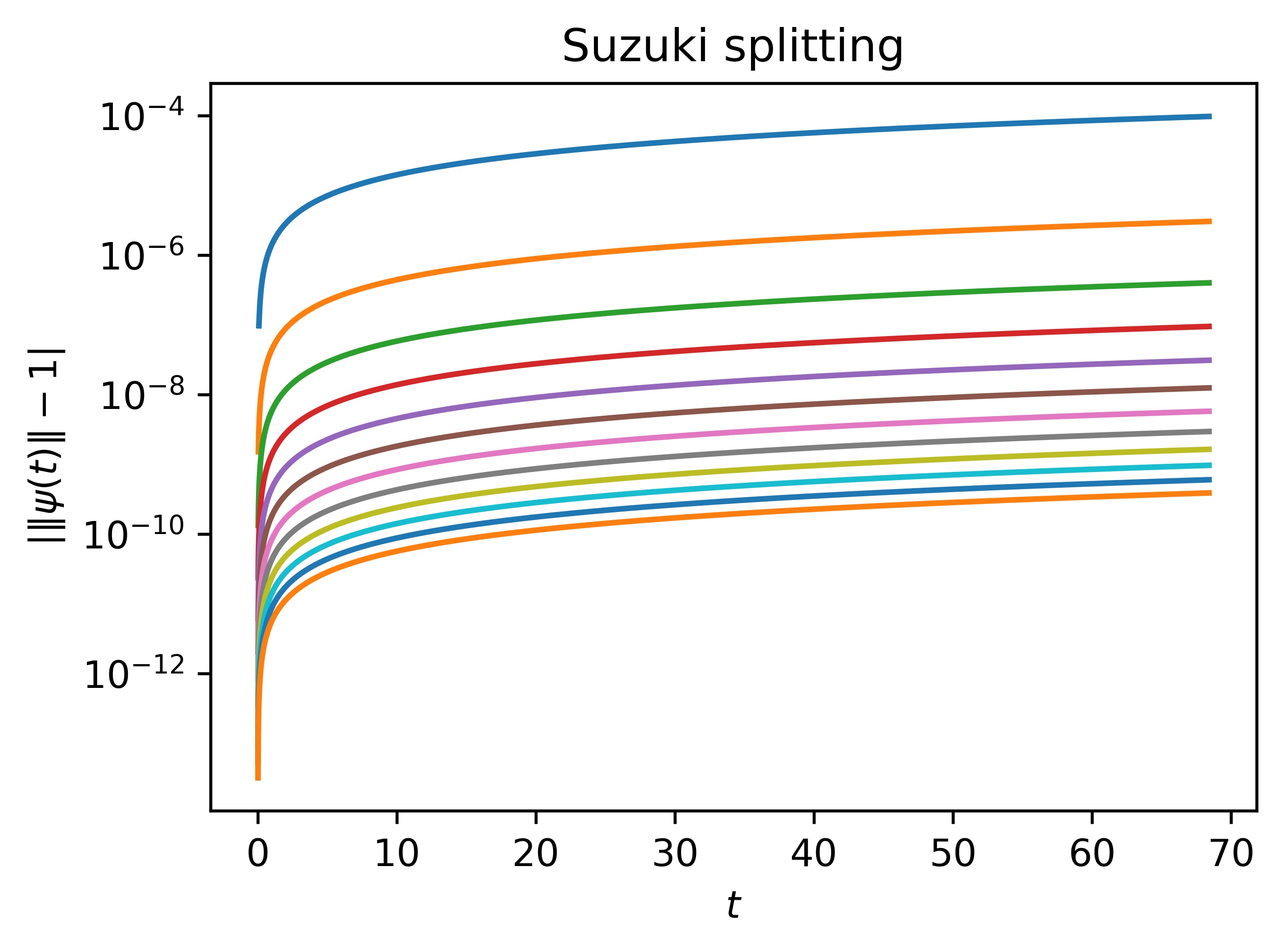} \\
\includegraphics[width=0.482\textwidth]{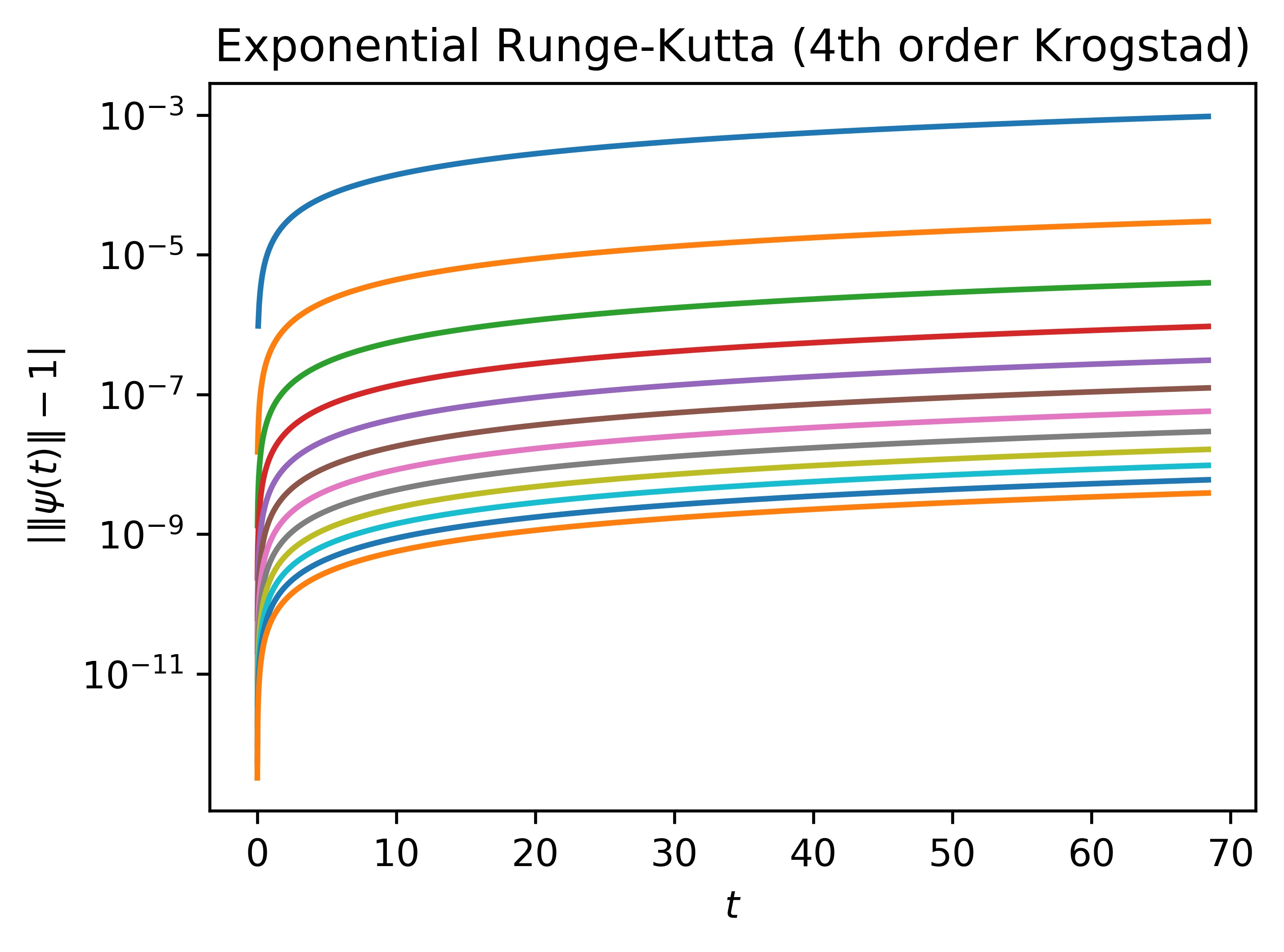}
\hspace{0.02\textwidth}
\includegraphics[width=0.482\textwidth]{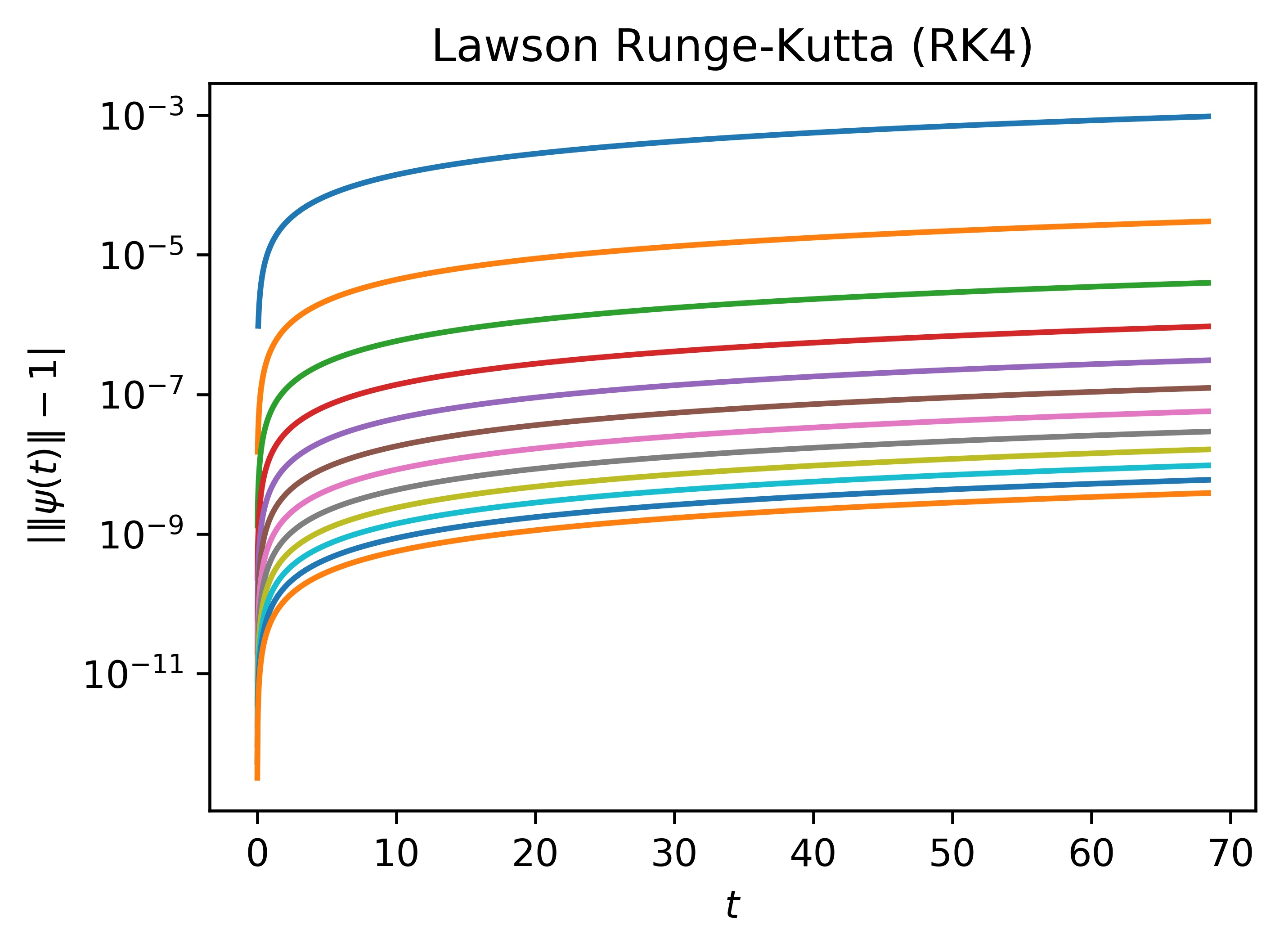} \\
\includegraphics[width=0.482\textwidth]{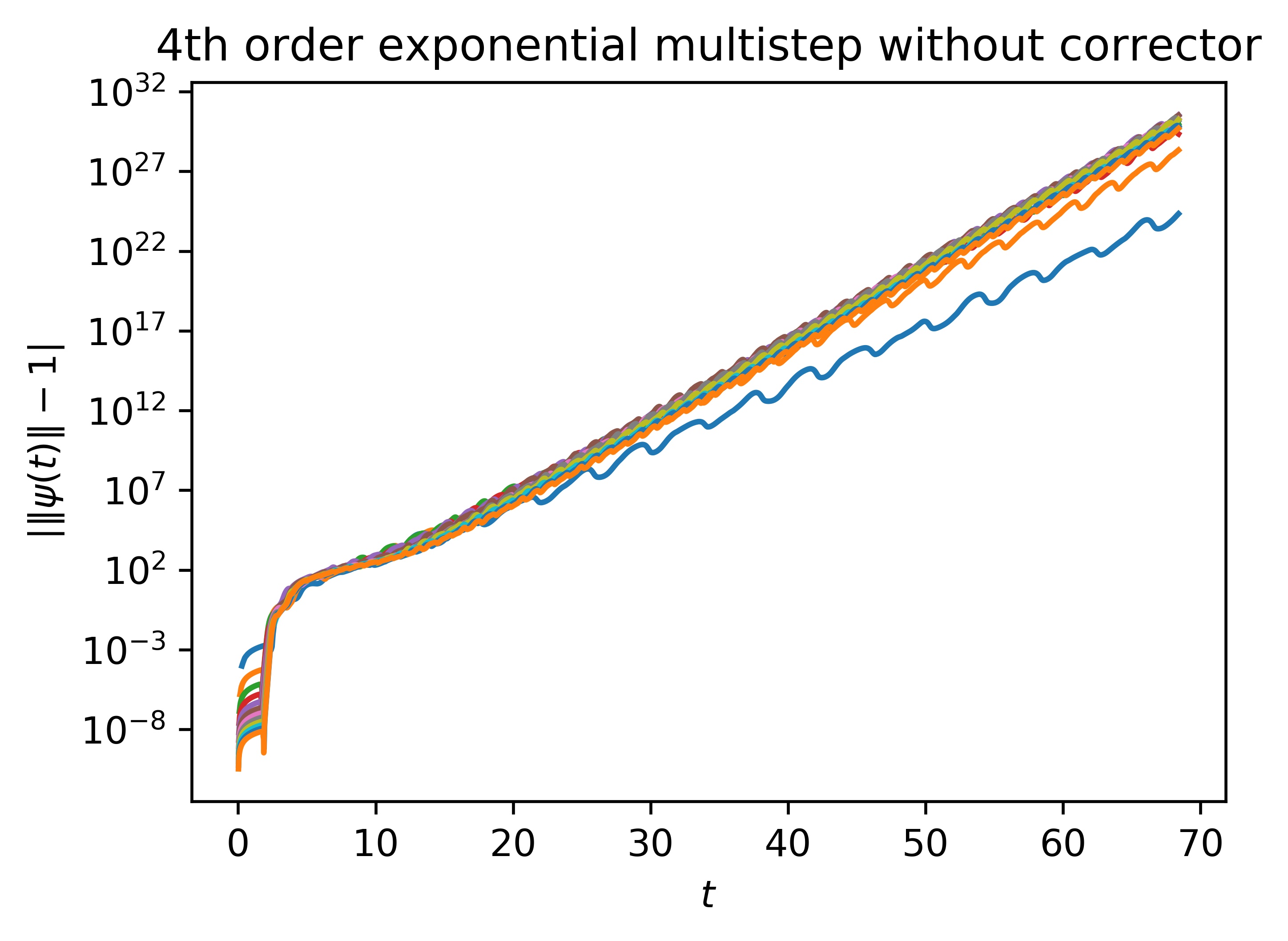}
\hspace{0.02\textwidth}
\includegraphics[width=0.482\textwidth]{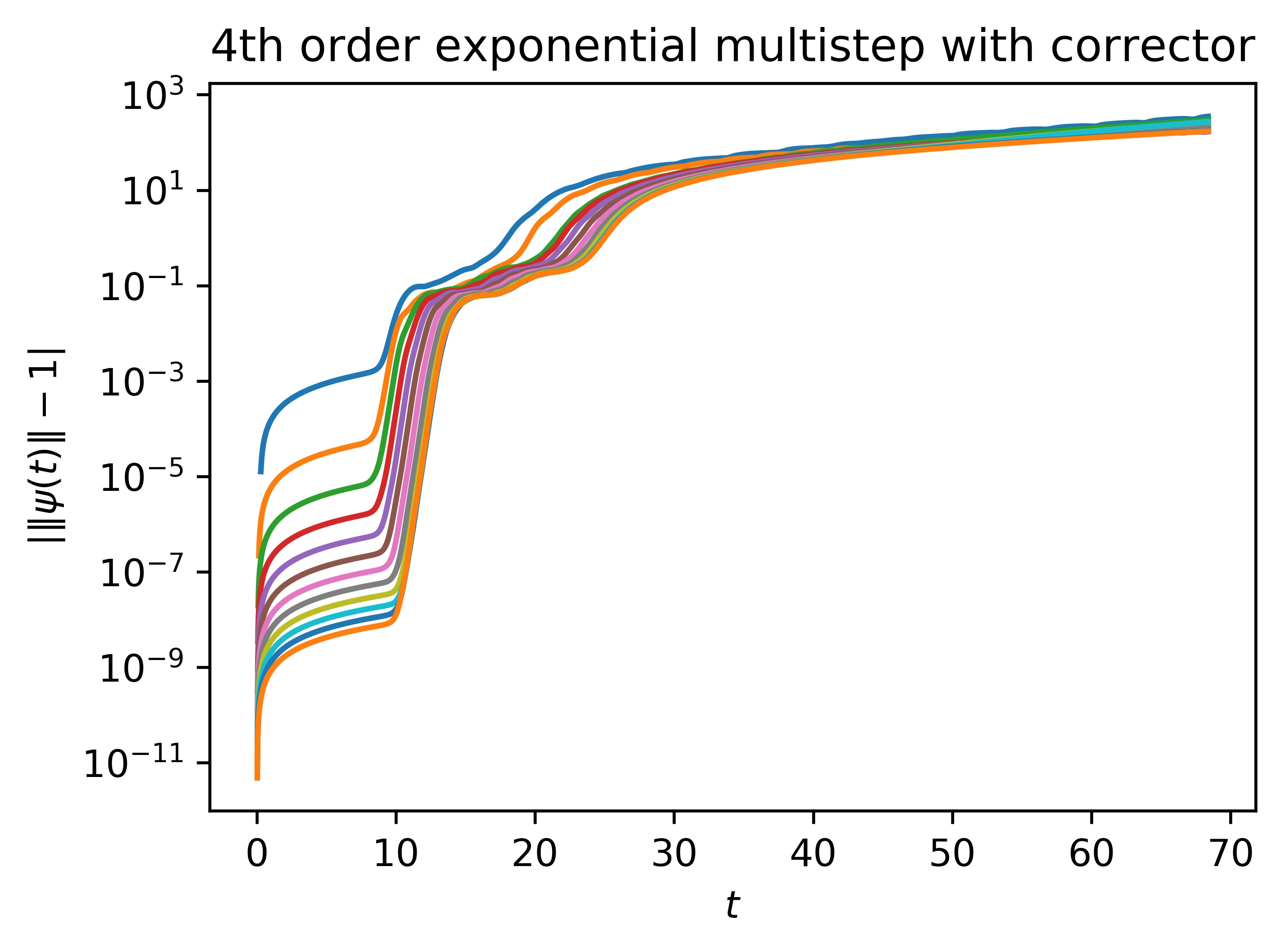} \\
\includegraphics[width=0.482\textwidth]{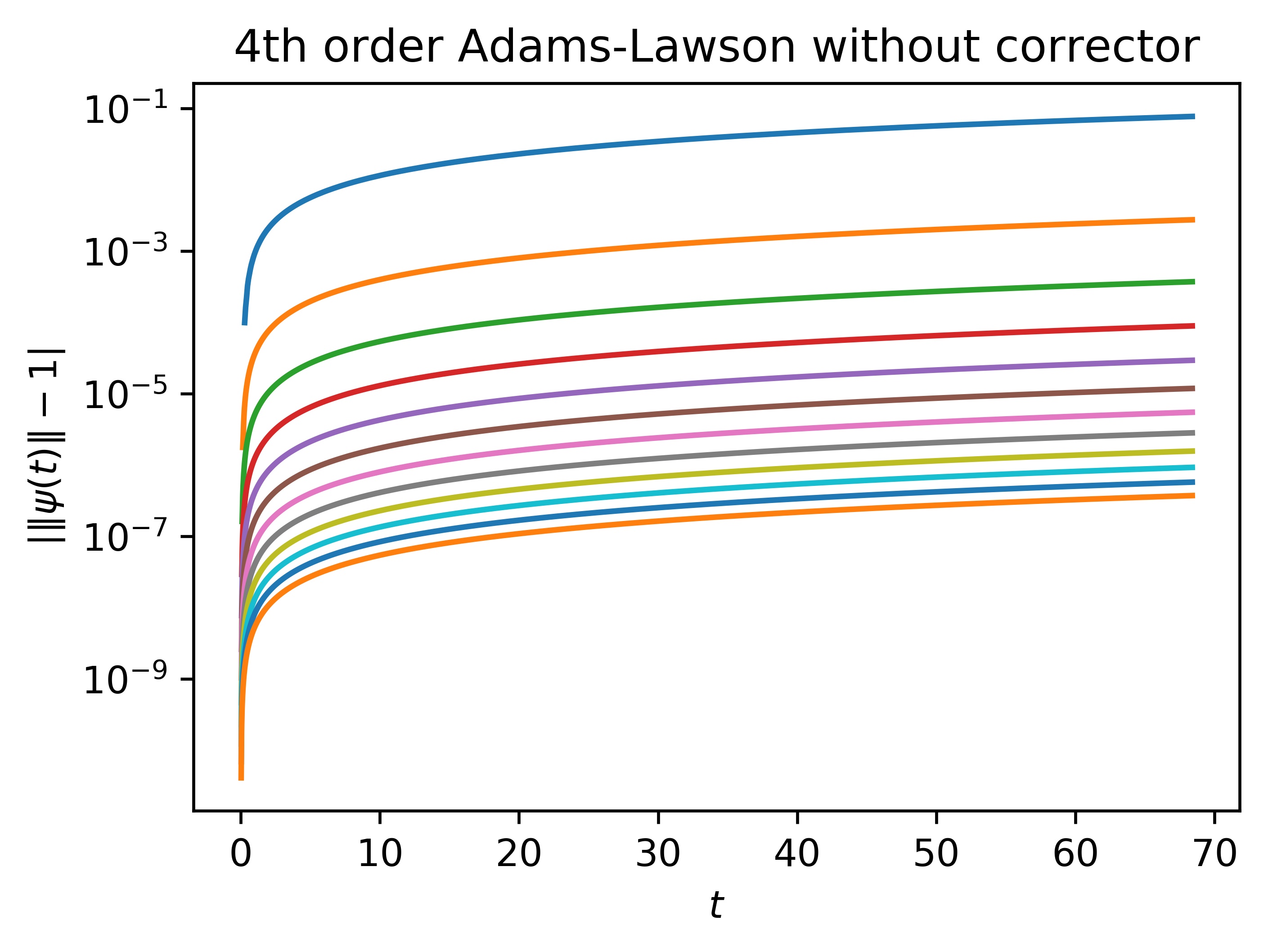}
\hspace{0.02\textwidth}
\includegraphics[width=0.482\textwidth]{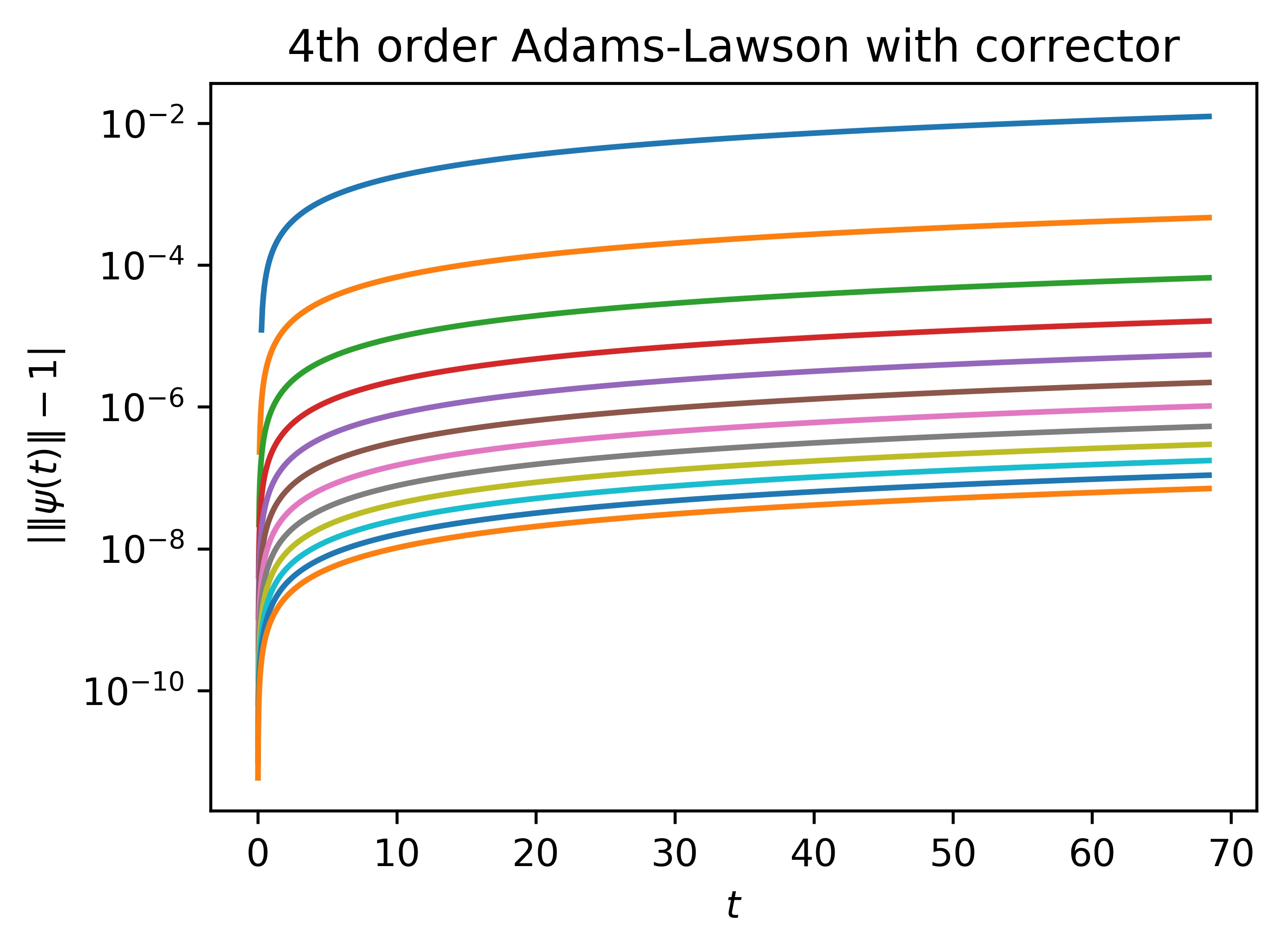} %\\
\caption{Comparison of stable long-time integration.\label{fig:stability}}
\end{center}
\end{figure}

To illustrate the performance of our numerical methods, we consider MCTDHF with the choice $N=4$
for a one-dimensional model of a helium atom \rev{investigated in~\cite{zanghellinietal03c}},
where\footnote{Note that in exponential
integrators, the explicit time-dependence in the potential does not call for a special treatment
in the numerical quadrature, in the splitting methods, the potential is propagated by
an explicit Runge-Kutta method of appropriate order.}
\begin{eqnarray*}
&&H(t) = H_0 + (x_1+x_2)\,\mathcal{E}(t), \\ %\label{irrad} \\
&&H_0  = -\frac{1}{2}(\partial_{x_1}^2+\partial_{x_2}^2)-\frac{2}{\sqrt{x_1^2+b^2}}
-\frac{2}{\sqrt{x_2^2+b^2}}+\frac{1}{\sqrt{(x_1-x_2)^2+b^2}},%\label{helium}
\end{eqnarray*}
with a smoothed Coulomb potential with shielding parameter $b=0.7408$, \rev{which is irradiated by a
short, intense, linearly polarized laser pulse}
$$\mathcal{E}(t) = \mathcal{E}_0\,g(t)\sin(\omega t).$$
The peak amplitude is set to $\mathcal{E}_0=0.1894$,
the frequency is $\omega=0.1837$, and we define the envelope
$g(t)=1.2\exp\big(-5\cdot 10^{-4}\left(t-6\pi/\omega\right)^2\big)$.
\rev{The parameters are taken from \cite{zanghellinietal03c}, and the envelope
is a smooth approximation of the trapezoidal envelope chosen there.
In \cite{zanghellinietal03c}, this model serves to illustrate the effect of
correlation on the probability density along the diagonal $x=y$,
which implies that the single-configuration Hartree--Fock approximation is insufficient.}
We first investigate stable long-time
propagation in Fig.~\ref{fig:stability}. We monitor norm conservation of the wave function
in the propagation of the ground state for the Hamiltonian $H_0=H(0)$ for
different equidistant stepsizes to resolve precisely the onset of instability.
For RK4, the number of steps is specified in the plot; for all other methods, the number of steps is in $\{1000,2000,\dots,12000\}$.
If norm conservation is violated beyond the effect of numerical accuracy, the method cannot be recommended
for physical applications. Indeed, we observe the following: Explicit Runge-Kutta methods only
behave in a stable way when the numerical accuracy is already very high, close to round-off error.
Exponential multistep methods\footnote{In this experiment, all multistep methods are started
by the Krogstad exponential Runge-Kutta method with stepsize $h/50$.} behave stably only for short times, even when a corrector step
is performed. Exponential Runge-Kutta and Runge-Kutta-Lawson methods behave
stably, likewise as splitting methods. Adams-Lawson multistep methods behave very stably,
a corrector step adds to the accuracy, as well as providing an error estimate as
the basis for adaptive time-stepping.
The unstable exponential multistep methods are no longer considered.
While showing the same stability behavior, the Yoshida splitting is demonstrated to be less efficient than
the Suzuki splitting, and the low order (but popular) Strang splitting is not competitive.
Higher-order multistep methods provide higher accuracy at the same computational effort irrespective of the order,
and are thus also considered for this comparison.

\begin{figure}[ht!]
\begin{center}
\includegraphics[width=11.5cm]{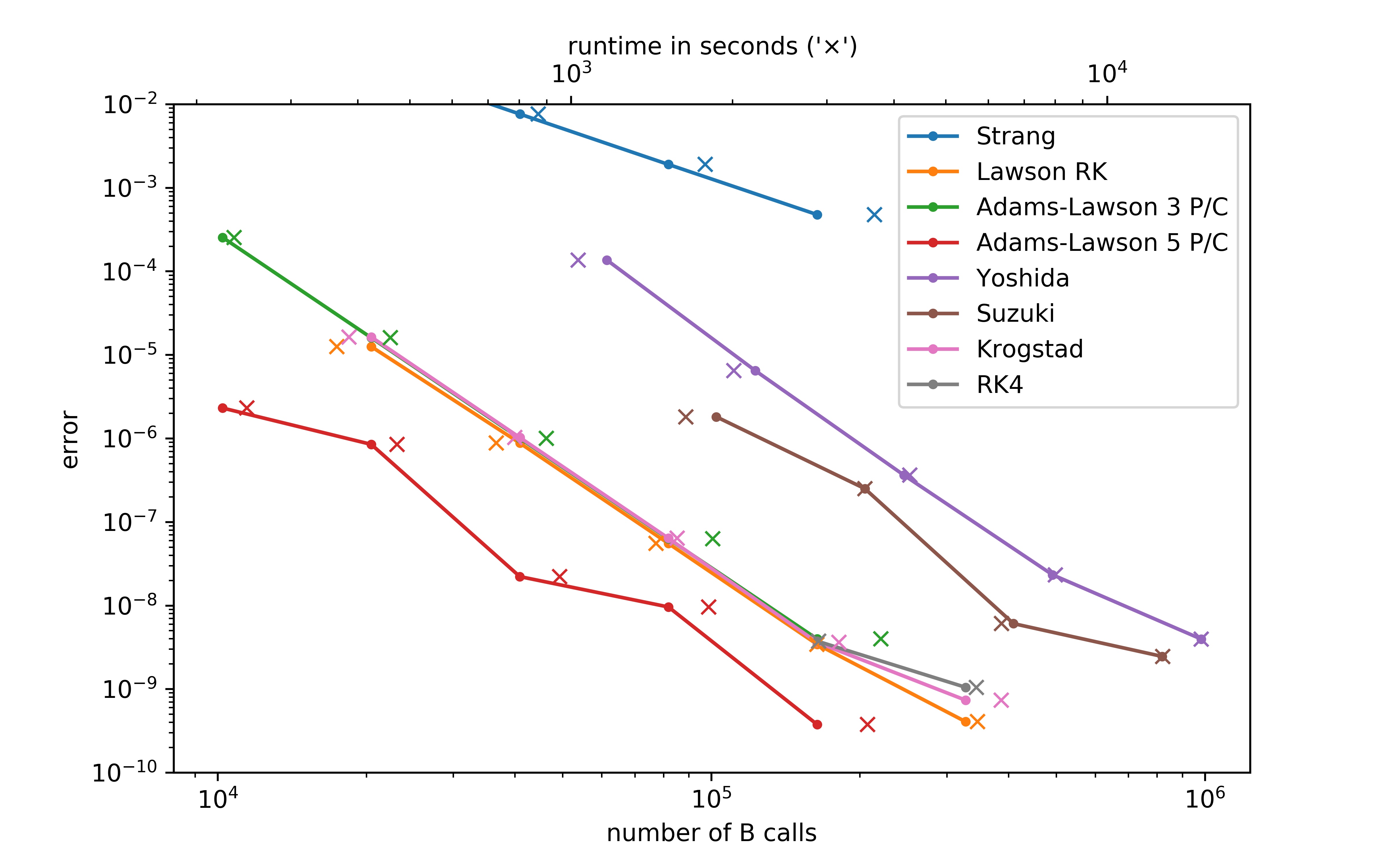}
\caption{Work/precision comparison for the helium atom, MCTDHF with $N=4$.
\rev{Dots on solid lines show the number of evaluations of $B$ and `$\times$' mark seconds of CPU time.}
\label{fig:bcalls}}
\end{center}
\end{figure}

Next, we compare the efficiency of the different integrators.
\rev{The unstable exponential multistep methods are no longer
considered. While showing the same stability behavior, the Yoshida splitting
is demonstrated to be less efficient than the Suzuki splitting, and the low order
(but popular) Strang splitting is not competitive. High-order multistep methods
provide higher accuracy at the same computational effort and are thus
also considered for this comparison.}
To this end, we plot
in Fig.~\ref{fig:bcalls} the accuracy as compared to a very precise reference solution
at $t=80$ as a function
of the number of evaluations of the computationally expensive potential part $B$ \rev{(dots on solid lines)}.
\rev{Furthermore, we give the CPU time required in a sequential implementation on one thread of the Vienna Scientific Cluster
(VSC) 3 comprising one Intel Xeon E5-2650v2 processor with 8 kernels of 2.6 gHz (crosses `$\times$').
We note that, as expected, the runtime is proportional to the number of potential evaluations.}
%, where in MCTDHF we have chosen $N=4$.
We observe that high-order Lawson \rev{multistep} methods perform best, where particularly the
high order which can be achieved in the multistep versions without additional evaluations
is advantageous.
Splitting methods, particularly the low order Strang splitting,
are not very efficient due to the high effort for the propagation of the potential.

\begin{figure}[hb!]
\begin{center}
\includegraphics[width=11.5cm]{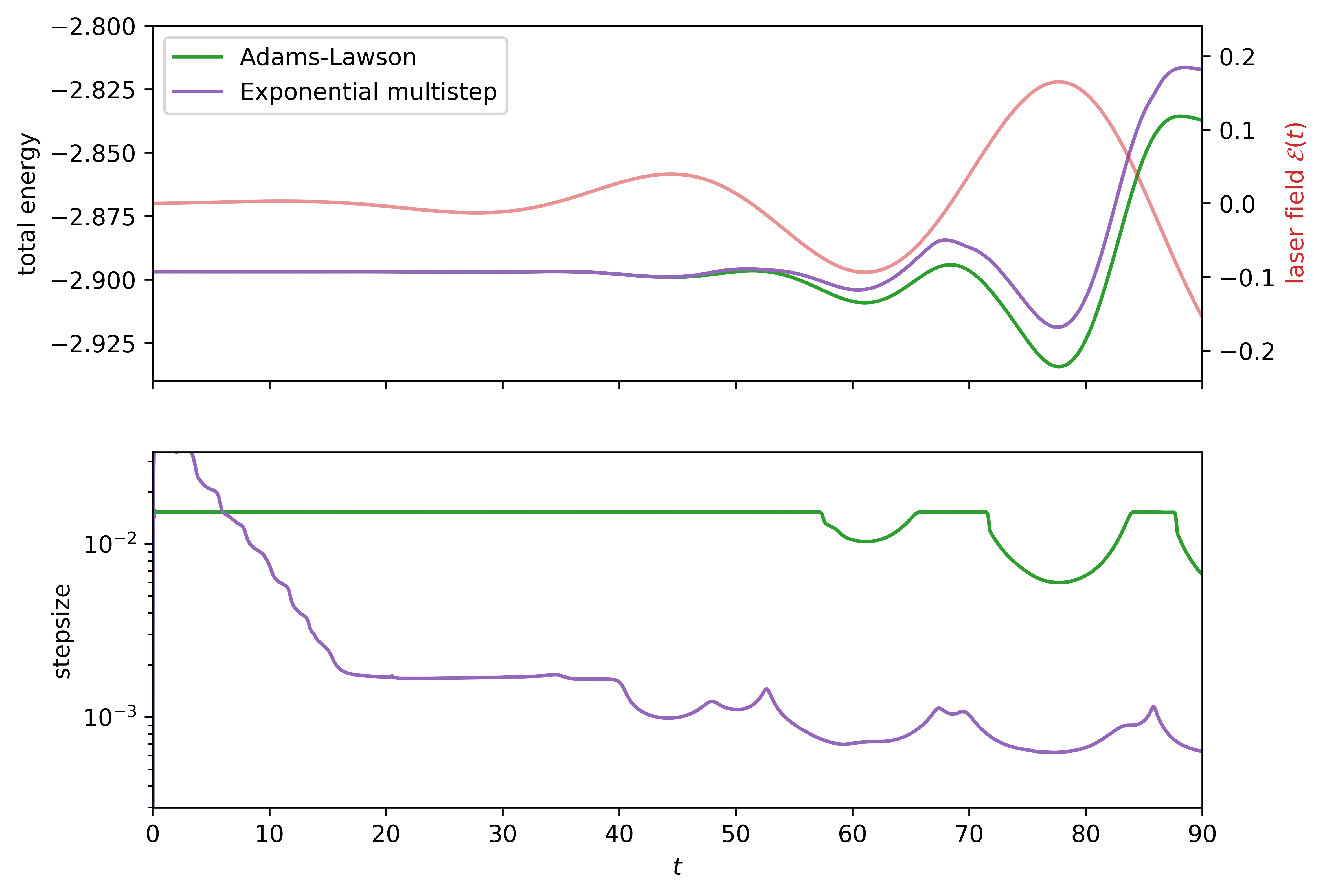}
\caption{Total energy functional and external potential (top) and automatically generated stepsizes (bottom),
MCTDHF with $N=4$.
\label{fig:stepsizes}}
\end{center}
\end{figure}

{We stress that this shows only the picture on uniform grids.
%Error estimators for the one-step methods have not been implemented due to the
%expected prohibitive computational effort for local error estimation.
The multistep versions
show their biggest advantage in adaptive time-stepping due to the cheap means of error estimation
in the predictor/corrector implementation.} To demonstrate that this
works reliably for Adams-Lawson methods, we show in Fig.~\ref{fig:stepsizes} the laser field $\mathcal{E}(t)$ and total energy functional (top)
illustrating the local solution smoothness, and the stepsizes (bottom)
automatically generated for the Adams-Lawson method of order 6.
We see that the adaptively chosen stepsizes reflect the smoothness of the time evolution
and the Lawson method enables larger stepsizes. The Adams-Lawson solution has been confirmed
to be converged to within the prescribed tolerance $10^{-5}$.
On the other hand, the corresponding exponential multistep method (without the Lawson transformation) shows
a noticeable deviation in the solution.

%\newpage

%\bibliographystyle{splncs04}
%\bibliography{num,schroedinger,books}

\begin{thebibliography}{10}
\providecommand{\url}[1]{\texttt{#1}}
\providecommand{\urlprefix}{URL }
\providecommand{\doi}[1]{https://doi.org/#1}

\bibitem{calvopal11}
Calvo, M., Palencia, C.: A class of explicit multistep exponential integrators
  for semilinear problems. Numer. Math.  \textbf{102},  367--381 (2011)

\rev{\bibitem{certaine60}
Certaine, J.: The solution of ordinary differential equations with large time
  constants. In: Ralston, A., Wilf, H. (eds.) Mathematical Methods for Digital
  Computers, pp. 128--132. Wiley, Hoboken, N.J. (1960)
}

\bibitem{frenkel34}
Frenkel, J.: Wave Mechanics, Advanced General Theory. Clarendon Press, Oxford
  (1934)

\bibitem{friedli78}
Friedli, A.: Verallgemeinerte {R}unge--{K}utta {V}erfahren zur {L}{\"o}sung
  steifer {D}ifferentialgleichungssysteme. In: Bulirsch, R., Grigorieff, G.,
  Schr{\"o}der, J. (eds.) Numerical Treatment of Differential Equations,
  Lecture Notes in Mathematics, vol.~631, pp. 35--50. Springer (1978)

\bibitem{haireretal02}
Hairer, E., Lubich, C., Wanner, G.: Geometric Numerical Integration.
  Springer-Verlag, Berlin--Heidelberg--New York (2002)

\bibitem{hocost10}
Hochbruck, M., Ostermann, A.: Exponential integrators. Acta Numer.
  \textbf{19},  209--286 (2010)

\bibitem{hocost17}
Hochbruck, M., Ostermann, A.: On the convergence of {L}awson methods for
  semilinear stiff problems. CRC Preprint 2017/9, KIT Karlsruhe Institute of
  Technology (2017),
  \texttt{https://www.waves.kit.edu/downloads/CRC1173\_Preprint\_2017-9.pdf}

\rev{\bibitem{koch19}
Koch, O.: Convergence of exponential {L}awson-multistep methods for the
  {MCTDHF} equations, to appear in M2AN Math. Model. Numer. Anal.
}

\bibitem{knth10a}
Koch, O., Neuhauser, C., Thalhammer, M.: Error analysis of high-order splitting
  methods for nonlinear evolutionary {S}chr\"{o}dinger equations and
  application to the {MCTDHF} equations in electron dynamics. M2AN
  Math.~Model.~Numer.~Anal.  \textbf{47},  1265--1284 (2013)

\bibitem{lawson67}
Lawson, J.: Generalized {R}unge--{K}utta processes for stable systems with
  large {L}ipschitz constants. SIAM J. Numer. Anal.  \textbf{4},  372--–380
  (1967)

\rev{\bibitem{noersett69}
N{\o}rsett, S.: An {A}-stable modification of the {A}dams--{B}ashforth methods.
  In: Conference on the Numerical Solution of Differential Equations, Lecture
  Notes in Mathematics, vol.~109, pp. 214–--219. Springer,
  Berlin--Heidelberg--New York (1969)
}

\bibitem{zanghellinietal03c}
Zanghellini, J., Kitzler, M., Brabec, T., Scrinzi, A.: Testing the
  multi-configuration time-dependent {H}artree-{F}ock method. J.~Phys.~B:
  At.~Mol.~Phys.  \textbf{37},  763--773 (2004)

\end{thebibliography}
\end{document}